\documentclass[11pt,a4paper]{article}
\usepackage[top=1in, bottom=1.25in, left=1in, right=1in]{geometry}
\usepackage[latin1]{inputenc}
\usepackage[english]{babel}
\usepackage{graphicx,epstopdf}
\usepackage{caption, subcaption,float}
\usepackage{color}
\usepackage{multirow}
\usepackage{tikz}
\usepackage{amsthm}
\usepackage{tabls}
\usepackage{cite}
\usepackage[ruled, linesnumberedhidden,vlined]{algorithm2e} 
\usetikzlibrary{matrix}
\usepackage{amsmath,amssymb,dsfont,mathtools}
\usepackage[square,sort,comma]{natbib}
\usepackage{hyperref}
\hypersetup{
    colorlinks=true,       
    linkcolor=blue,          
    citecolor=olive,       
    filecolor=magenta,      
    urlcolor=blue          
}

\begin{document}

\title{On the loxodromic actions of Artin-Tits groups}
\date{\today }
\author{Mar\'{i}a Cumplido\\  Universit\'{e} Rennes 1, Universidad de Sevilla}

\maketitle

\newtheorem{theorem}{Theorem}
\newtheorem{lemma}[theorem]{Lemma}
\newtheorem{proposition}[theorem]{Proposition}
\newtheorem{corollary}[theorem]{Corollary}

\theoremstyle{definition}
\newtheorem{definition}[theorem]{Definition}

\theoremstyle{remark}
\newtheorem{remark}[theorem]{Remark}

\newcommand{\St}[1]{\operatorname{St}(\mathcal{#1})}
\newcommand{\po}{\preccurlyeq}
\newcommand{\wedger}{\wedge^\Lsh}
\newcommand{\veer}{\vee^\Lsh}
\newcommand{\bigveer}{\bigvee^\Lsh}
\newcommand{\alphae}{\alpha_{\text{ext}}}
\newcommand{\gammae}{\gamma_{\text{ext}}}
\newcommand{\notpo}{\npreceq}
\newcommand{\so}{\succcurlyeq}
\newcommand{\notso}{\nsucceq}
\newcommand{\lox}[1]{\operatorname{Lox}(#1)}
\newcommand{\ca}{\mathcal{C}_{AL}}

\begin{abstract}
Artin-Tits groups act on a certain delta-hyperbolic complex, called the ``additional length complex". For an element of the group, acting \emph{loxodromically} on this complex is a property analogous to the property of being pseudo-Anosov for elements of mapping class groups. By analogy with a well-known conjecture about mapping class groups, we conjecture that ``most'' elements of Artin-Tits groups act loxodromically. More precisely, in the Cayley graph of a subgroup $G$ of an Artin-Tits group, the proportion of loxodromically acting elements in a ball of large radius should tend to one as the radius tends to infinity. In this paper, we give a condition guaranteeing that this proportion stays away from zero. This condition is satisfied e.g. for Artin-Tits groups
of spherical type, their pure subgroups and some of their commutator subgroups.
\end{abstract}

\section{Introduction}

Let $A$ be an irreducible Artin-Tits group of spherical type which acts on the additional length graph, denoted $\ca$, which is a $\delta$-hyperbolic complex introduced in \citep{Calvez2016graph}. The Cayley graph of a subgroup $G$ of $A$ with generator system $S$ will be denoted $\Gamma(G,S)$.  The ball on $\Gamma(G,S)$ centered in the trivial vertex with radius $R$ will be denoted by $B_{\Gamma(G,S)}(R)$, and by $\lox{G,\ca}$ we mean the set of all the elements in $G$ that act loxodromically on $\ca$. 

\medskip

There is a well-known conjecture which claims that ``most'', or ``generic'' elements of the mapping classes of a surface are pseudo-Anosov: picking an element of the mapping class group ``at random'' yields a pseudo-Anosov element with overwhelming probability. The braid group with $n$ strands happens to be both a mapping class group (of the $n$-punctured disk) and the Artin-Tits group $A_{n-1}$, and it was proven in \citep{Calvez2016graph} that if a braid acts loxodromically on $\ca$, then it is pseudo-Anosov (this is actually believed to be an equivalence). By analogy, this justifies the following conjecture: ``generic'' elements of an Artin-Tits group $A$ (or of a reasonable subgroup $G$) act loxodromically on $\ca$.

\smallskip

In order to make the definition of genericity more precise, let us describe two possible definitions of picking randomly an element in $\Gamma(G,S)$. The first method consists of performing a random walk of a certain length~$\ell$ on~$\Gamma(G,S)$, starting from the identity vertex; the second method is to pick a random vertex from the ball $B_{\Gamma(G,S)}(R)$ of radius~$R$. 

\smallskip

The \emph{genericity conjecture} claims that, if~$G$ contains at least one element which acts loxodromically on~$\mathcal C_{AL}$, then picking randomly an element in $\Gamma(G,S)$ yields a loxodromically acting element with a probability that tends to one when the length~$\ell$ of the random walk, or when the radius~$R$ of the ball tends to infinity. In the ``large balls" model, this is claiming that
$$
\lim_{R\to\infty} \dfrac{|\lox{G,\ca} \cap B_{\Gamma(G,S)}(R) |}{|B_{\Gamma(G,S)}(R)|} \ =1.
$$
The paper \citep{Calvez2016} contained a partial proof of this conjecture, namely when $G~=~A$, and
\begin{itemize}
\item either genericity was defined according to the random walks model,
\item or genericity was defined according to the large balls model, but only with respect to one very particular generating set~$S$, namely Garside's set of simple elements (see below).
\end{itemize}

\medskip

In the present paper (Lemma~\ref{tecnique}), we will give a condition for $G$  so that the proportion of loxodromic elements in $B_{\Gamma(G,S)}(R)$ stays away from zero when $R$ tends to infinity. Notice that this is a weaker conclusion, but it will be proven for \emph{every} generator system $S$. Our main result is the following:

\begin{theorem}\label{main}

Let $G$ be an Artin-Tits group of spherical type, a pure subgroup of an Artin-Tits group of spherical type or the commutator subgroup of $I_{2(2m+1)}$, $A_n,\, D_n,\, E_n$ or $H_n$. Let $S$ be any generator set of $G$. Then the following condition is satisfied:

$$\liminf_{R\to\infty} \dfrac{|\lox{G,\ca} \cap B_{\Gamma(G,S)}(R) |}{|B_{\Gamma(G,S)}(R)|} \ >0$$
 
\end{theorem}

\section{Reminders}

In this section we will recall the main concepts and results used in this paper.  First, we define the groups that participate in Theorem~\ref{main}.

\begin{definition}[Artin-Tits group] \label{artin}
Let $I$ be a finite set and $M=(m_{i,j})_{i,j\in I}$ a symmetric matrix with $m_{i,i}=1$ and $m_{i,j}\in\{2,\dots, \infty \}$ for $i\neq j$. Let $S=\{\sigma_i\,|\, i\in I\}$. The Artin-Tits system associated to $M$ is $(A,S)$, where $A$ is a group with the following presentation 
$$A=\langle S \,|\, \underbrace{\sigma_i\sigma_j\sigma_i\dots}_{m_{i,j} \text{ elements}}=\underbrace{\sigma_j\sigma_i\sigma_j\dots}_{m_{i,j} \text{ elements}} \forall i,j\in I,\, i\neq j,\, m_{i,j}\neq \infty \rangle.$$
\end{definition}

\noindent The Coxeter group $W$ associated to $(A,S)$ can be obtained by adding the relations $\sigma_i^2=1$:
$$W_A=\langle S \,|\, \sigma_i^2=1 \, \forall i\in I ; \underbrace{\sigma_i\sigma_j\sigma_i\dots}_{m_{i,j} \text{ elements}}=\underbrace{\sigma_j\sigma_i\sigma_j\dots}_{m_{i,j} \text{ elements}} \forall i,j\in I,\, i\neq j,\, m_{i,j}\neq \infty \rangle.$$
An Artin-Tits system $(A,S)$ can be represented with a \emph{Coxeter graph}, denoted~$\Gamma_A$. The set of vertices of $\Gamma_A$ is $S$, and there is an edge joining two vertices $s,t\in \Sigma$ if $m_{s,t}\geq 3$. The edge will be labelled with $m_{s,t}$ if $m_{s,t}\geq 4$.  

\smallskip

\begin{figure}[h]
  \centering
  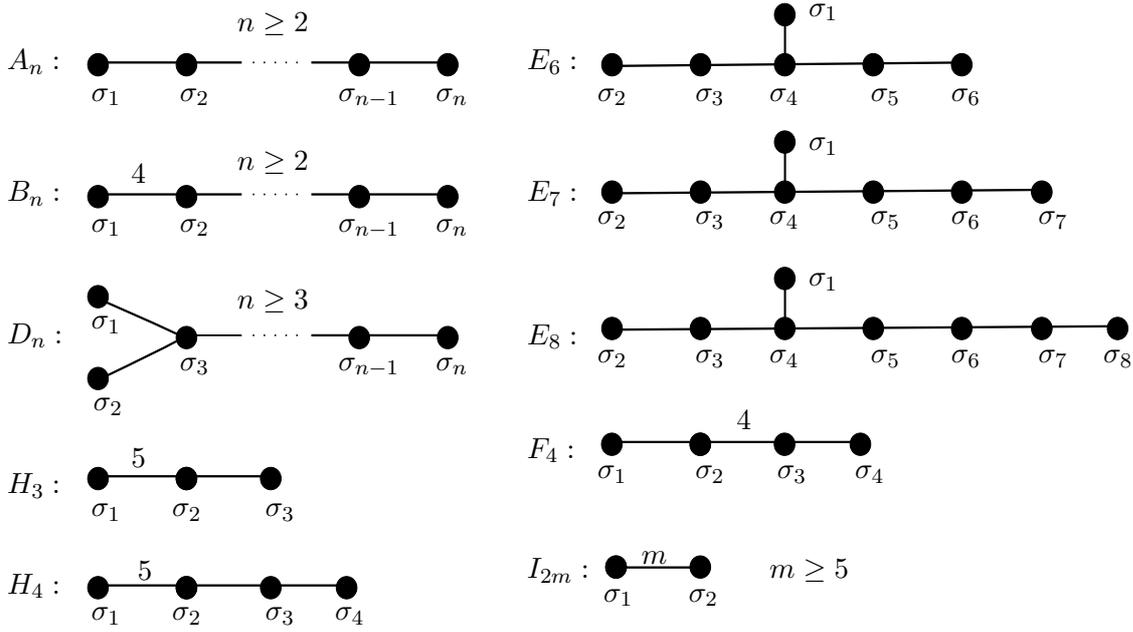
  \medskip
  \caption{Classification of Artin-Tits groups of spherical type.}
  \label{coxetergraphs}
\end{figure}

Through this paper, we shall only be interested in Artin-Tits groups of spherical type, meaning that their associated Coxeter groups are finite. It is well known \citep{Coxeter} that Artin-Tits group of spherical type are classified into 10 different types. The classification is described in Figure~\ref{coxetergraphs} by using Coxeter graphs.

\bigskip

The best-known example for these groups is the \emph{braid group} with $n$ strands ($A_{n-1}$ in Figure~\ref{coxetergraphs}). Each braid with $n$ strands can also be seen as a collection of $n$ disjoint paths in a cylinder, defined up to isotopy, joining $n$ points at the top with $n$ points at the bottom, running monotonically in the vertical direction. In this case, each generator $\sigma_i$ represents a crossing between the strands in the position $i$ and $i+1$ with a fixed orientation, while $\sigma_i^{-1}$ represents the same crossing with the opposite orientation.

\begin{definition}
Let $(A,S)$ be an Artin-Tits system where $A$ is an Artin-Tits group of spherical type. Its \emph{pure} subgroup $P(A)\subset A$ is defined as follows

$$P(A):= \ker(p: A \longrightarrow W_A) $$ where $p$ is the canonical projection from $A$ to its associated Coxeter group $W_A$.

\end{definition}

The techniques used in the proofs that we will see involve Garside theory.  Let us now recall some important concepts about this theory. A group $G$ is called a \emph{Garside group} with Garside structure $(G,\mathcal{P},\Delta)$ if it admits a submonoid $\mathcal{P}$ of positive elements such that $\mathcal{P}\cap \mathcal{P}^{-1}=\{1\}$ and a special element $\Delta \in \mathcal{P}$, called Garside element, satisfying the following:

\begin{itemize}
 \item There is a partial order in $G$, $\po$,  defined by $a \po b \Leftrightarrow a^{-1}b \in \mathcal{P}$ such that for all $a,b\in G$ there exists a unique gcd $a \wedge b$ and an unique lcm $a \vee b$.  This order is called prefix order and it is invariant by left-multiplication.

\item The set of simple elements $[1,\Delta]=\{a\in G\,|\, 1\po a \po \Delta  \}$ generates G.

\item $\Delta^{-1}\mathcal{P} \Delta = \mathcal{P}$. 

\end{itemize}

In a Garside group, the monoid $\mathcal{P}$ also induces a partial order invariant under \emph{right}-multiplication. This is the suffix order, $\so$, defined by $a\so b \Leftrightarrow ab^{-1}\in \mathcal{P}$, such that for all $a,b\in G$ there exists an unique lcm $(a\wedger b)$ and an unique gcd $(a\veer b)$. It is well known that every Artin-Tits group of spherical type admits a Garside structure of \emph{finite type} \citep{Brieskorn1972, Dehornoy1999}, which means that $[1,\Delta]$ is finite.

An \emph{atom} is an element $a\in \mathcal{P}$ such that for every $b\in\mathcal{P}$, $b\not\po a$. The atoms of an Artin-Tits group of spherical type, $A$, equipped with the Garside structure mentioned above, are the generators of the presentation given in Definition~\ref{artin}. 

\begin{definition}[Normal forms]
Given two simple elements $a,b$, the product $a\cdot b$ is said to be in left (respectively right) normal form if $ab\wedge \Delta =a$ (respectively $ab\wedger \Delta =b$). 

We say that $ \Delta^k x_1\cdots x_r$ is the \emph{left normal form} of an element $x$ if $k\in \mathbb{Z}$, $x_i\notin \{1,\Delta\}$ is a simple element for $i=1,\ldots , r$, and $x_i x_{i+1}$ is in left normal form for $0<i < r$. Analogously, $x_1\cdots x_r\Delta^k$ is the \emph{right normal form} of $x$ if $k\in \mathbb{Z}$, $x_i\notin \{1,\Delta\}$ is a simple element for $i=1,\ldots , r$, and $x_i x_{i+1}$ is in right normal form for $0<i < r$. 

It is well known that the normal forms of an element are uniques \citep{Dehornoy1999} and that the numbers $r$ and $k$ do not depend on the normal form (left or right). We define the infimum, the canonical length and the supremum of $x$ respectively as  $\inf(x)~=~k$, $\ell(x)=r$ and $\sup(x)=k+r$.

\end{definition}

\begin{remark}
Notice that, even though we will work with an arbitrary generator system $S$ for a subgroup $G$ of $A$, every time we talk about normal forms we will be always using the above definition, where each letter $x_i$ of the normal form is simple and maybe not in $S$. 
\end{remark}

\begin{definition}[Rigidity] Let  $x = \Delta^k x_1\cdots x_r$ be in left normal form. We define the initial and the final factor respectively as $\iota(x)=\Delta^{k} x_1 \Delta^{-k}$ and $\varphi(x)=x_r$. We will say that $x$ is \emph{rigid} if $\varphi(x)\cdot \iota(x)$ is in left normal form. 

\end{definition}

\section{Proportion of loxodromic actions}

In this section we will give the condition to have the positive proportion of loxodromic actions of Artin-Tits group of spherical type and we will apply it to some subgroups. From now on, suppose that $A$ is an Artin-Tits group of spherical type.

\begin{lemma}[{\citealp{Calvez2016}}]\label{xg} 
For every atom $a\in A$, there is an element $x_a\in A$ which acts loxodromically on $\ca$ such that the left and the right normal forms of $x_a$ are the same and $\iota(x_a)=\varphi(x_a)=a$.

Moreover, if $g\in A$ is rigid and its normal form contains the subword $w_a:=x_a^{390}$, then $g$ acts loxodromically on $\ca$.
\end{lemma}

The previous lemma is a summary of some results proven by Calvez and Wiest. For every atom $a$, the element $x_a$ is constructed in \citep[Theorem 3.1]{Calvez2016}. The second part of the lemma is proven in \citep[Lemma 5.3]{Calvez2016}.

\begin{lemma}\label{tecnique}
Let $G$ be a subgroup of $A$. Suppose that there is a finite set $X$ of elements in $G$ such that for every $g\in G$ there exists $x\in X$ such that $g\cdot x$ is rigid and its normal form contains the subword $w_a$. Let $B(R)~:=~B_{\Gamma(G,S)}(1,R)$, where $S$ is a finite generator system of $G$. Then there are constants $\varepsilon , R_0 > 0$ depending on $S$, such that for all $R > R_0$,

$$\dfrac{|\lox{G,\ca} \cap B (R) |}{|B (R)|} > \varepsilon $$

\end{lemma}

\emph{Proof.} Let $R_0=\max\{\ell_{S}(x)\,|\, x\in X\}$ be the maximum of the canonical lengths, with respect to the generator system $S$, of all the elements in~$X$. By Lemma~\ref{xg}, $x\cdot g$ acts loxodromically on $\ca$, for some $x\in X$. Thus, $d(g,\lox{G,\ca})\leq R_0$, for every~$g\in~G$. In particular for every~$g \in B(R-R_0)$, there is a loxodromic element which is at distance at most~$R_0$ from~$g$ and which lies in~$B(R)$. Therefore,  $$|B(R-R_0)|\leq |\lox{G, \ca} \cap B(R)|\cdot |B(R_0)|,$$

which implies that 

$$\dfrac{|\lox{G, \ca} \cap B(R)|}{|B(R)|}\geq \dfrac{|B(R-R_0)|}{|B(R)| }\cdot\dfrac{1}{B(R_0)}.$$
%Hence,
%
%
%$$\dfrac{|\lox{G,\ca} \cap B(1,R)|}{|B(1,R-R_0)|}\geq \frac{1}{|B(1,R_0)|}.$$

Moreover, since the number of elements of a ball in a Cayley graph grows at most exponentially with its radius, there is a constant $\varepsilon'$ which depends on $S$ and $R_0$ such that 

$$\dfrac{|B(R-R_0)|}{|B(R)|}\geq \varepsilon',$$
for all $R$. Thus, the number $\varepsilon= \dfrac{\varepsilon'}{|B(R_0)|}$ does the job.
\hfill $\blacksquare$

\medskip

\subsection{Proof of Theorem~\ref{main}}

Theorem~\ref{main} is a summary of the following three theorems, which will be proven separately. 

\begin{theorem}
Let $(A,S)$ be an Artin-Tits system and let $S'$ be a generator system for $A$. Let $B(R):=B_{\Gamma(A,S')}(1,R)$. Then there are constants $\varepsilon , R_0 > 0$ depending on~$S'$, such that for all $R > R_0$,

$$\dfrac{|\lox{A,\ca} \cap B (R) |}{|B (R)|} > \varepsilon $$

\end{theorem}

\emph{Proof.} This is now an immediate consequence of Lemma~\ref{xg} and Lemma~\ref{tecnique}, applied with~$G~=~A$. \hfill $\blacksquare$

We want to prove this also for pure subgroups of Artin-Tits groups. In particular, the following theorem shows that the pure braid group has a positive proportion of pseudo-Anosov elements.

\begin{theorem}\label{pure}
Let $G\subseteq A$ be the pure subgroup of an Artin-Tits group, equipped with any finite generator system $S$. Define $B(R):=B_{\Gamma(G,S)}(1,R)$ . Then there are constants $\varepsilon , R_0 > 0$ depending on $S$, such that for all $R > R_0$,

$$\dfrac{|\lox{G,\ca} \cap B (R) |}{|B (R)|} > \varepsilon $$

\end{theorem}

%Notice that $S$ is finite because $G$ is of finite index on $A$.

\emph{Proof.} We want to prove that Lemma~\ref{tecnique} can be applied.  We claim that for each $g\in G$ we can find an $x\in G$ such that $g\cdot x$ is rigid and its left normal form contains the subword $w_a$. The choice of the atom $a$ depends on $A$:  we choose specifically $a=\sigma_2$ for $B_n$, $H_3$, $H_4$, $F_4$, $I_{2m}$, $a=\sigma_3$ for $D_n$ and $a=\sigma_4$ for $E_6$, $E_7$, $E_8$. For $A_n$ we can set $a$ equal to any atom $\sigma_i$, $i=1,\dots, n-1$. We also have to prove that the length of~$x$ is bounded from above, in order to guaranteed the finiteness of the family~$X$ described in Lemma~\ref{tecnique}. The desired element $x$ will be constructed as a product $x=w_z\cdot w_a\cdot w_p\cdot w_r$, where $w_z$, $w_a$, $w_p$ and $w_r$ are words in left normal form with infimum equal to zero. We also need the whole word $w_g\cdot w_z\cdot w_a\cdot w_p\cdot w_r$ to be in left normal form, where $w_g$ is the the left normal form of $g$. Before going into details, we describe the function of each factor of $x$. The first word $w_z$ makes sure that the product $g\cdot x$ is in left normal form; $w_a$ is the element mentioned in Lemma~\ref{tecnique}; $w_p$ provides the pureness of $x$; finally, $w_r$ provides the rigidity of the element $gx$. 

\medskip

We recall that by Lemma~\ref{xg}, $\iota(w_a)=\phi(w_a)=a$. Let $b$ be an atom such that $\phi(w_g)\so b$. Let $p: A\longrightarrow W_A$  be the canonical projection of $A$ to its Coxeter group and let $s$ be an atom such that $s\not\po \iota(g)$. Now, the aim is to write
$$w_z= \boldsymbol{b}\cdot w_z'\cdot \boldsymbol{a},\quad w_p= \boldsymbol{a}\cdot w_p'\cdot \boldsymbol{a},\quad w_r= \boldsymbol{a}\cdot w_r'\cdot \boldsymbol{\Delta s^{-1}}.$$ 
where the first and the last factor of the normal forms of $w_z$, $w_p$ and $w_r$ are distinguished in bold. In order to help visualize the idea of the proof, we say that $x$ should be such that the normal form of $g\cdot x$ is

\[ \underbrace{\overbrace{\boldsymbol{\iota(g)}}^{s\not\po}\cdots \overbrace{\boldsymbol{\varphi(g)}}^{\so b}}_{w_g} \cdot \underbrace{\boldsymbol{b}\cdots \boldsymbol{a}}_{w_z} \cdot \underbrace{\boldsymbol{a}\cdots \boldsymbol{a}}_{w_a} \cdot \underbrace{\boldsymbol{a}\cdots \boldsymbol{a}}_{w_p} \cdot \underbrace{\boldsymbol{a}\cdots (\boldsymbol{\Delta s^{-1}})}_{w_r}\]

Words $w_z$ and $w_r$ with these characteristics can be constructed with less than 6 factors. This is proven in \citep[Lemma 3.4]{Caruso} in the case of braids and in \citep[Propositions 57-65]{Gebhardt2016} for the other Artin-Tits group of spherical type.

\medskip
The last step is to construct, for any given $w\in W_A$, a word in left normal form  $w_p$ whose first and last factors are equal to $a$ and which is such that $p(w_p)=w$. So, if we consider a generator system $S'$ for $W_A$, then we need to find for each $s'\in S'$ an element $\tilde{s}\in p^{-1}(s')$ such that the first and last factor of its formal form is $a$. Hence, with the elements of the form $\tilde{s}$ it is possible to construct $w_p$ with the desired properties. In order to construct this element $\tilde{s}$, let $$\Sigma_{i,j}=\left\{\begin{array}{ll}
\sigma_i\sigma_{i+1}\cdots\sigma_j & \text{if } i<j \\
\sigma_i\sigma_{i-1}\cdots\sigma_j & \text{if } j<i 
\end{array}  \right. , \quad
\Sigma_{i,j}^{(2)}=\left\{\begin{array}{ll}
\sigma_i^2\sigma_{i+1}^2\cdots\sigma_j^2 & \text{if } i<j \\
\sigma_i^2\sigma_{i-1}^2\cdots\sigma_j^2 & \text{if } j<i 
\end{array}  \right.$$

We study in Table~\ref{tablaprueba} each possible class for $A$ (see Figure~\ref{coxetergraphs}). Some remarks about this table are given below.

\begin{table}[]
\begin{tabular}{|c||c|c|c|c|c|c||c|c|c|}
\hline
$A$         & \multicolumn{6}{c||}{$\boldsymbol{A_n}$} & \multicolumn{3}{c|}{$\boldsymbol{B_n}$}                                                                                                                   \\ \hline
$a$         & \multicolumn{2}{c|}{$\sigma_1$}                         & \multicolumn{2}{c|}{$\sigma_i,\,\, i=2,\dots n-1$}                            & \multicolumn{2}{c||}{$\sigma_n$}       &            \multicolumn{3}{c|}{$\sigma_2$}                                                                                                                   \\ \hline 
$s'$         &  $\sigma_1$ & $\Sigma_{n,1}$                               & $\sigma_i$ & $\Sigma_{1,n}$                                                 & $\sigma_n$ & $\Sigma_{1,n}$   & $\sigma_1$                     & $\sigma_2$ & $\Sigma_{2,n}$                       \\ \hline  
$\tilde{s}$ & $\sigma_1$ & $\,\,\Sigma^{(2)}_{1,n}\Sigma_{n,1}\sigma_1^2$ & $\sigma_i$ & $\,\,\Sigma_{i,1}^{(2)}\Sigma_{1,n}\Sigma_{n,i}^{(2)}$ & $\sigma_n$ & $\,\Sigma_{n,1}^{(2)}\Sigma_{n,1}\sigma_1^2\,$ & $\,\sigma_2^2\sigma_1\sigma_2^2\,$  &  $\sigma_2$ & $\,\sigma_2^2\Sigma_{2,n}\Sigma^{(2)}_{n,2}\,$ \\ \hline
\end{tabular}

\medskip

\begin{tabular}{|c||c|c|c||c|c|c||c|c|c|}
\hline
$A$       & \multicolumn{3}{c||}{$\boldsymbol{D_n}$}                                                                                                                      & \multicolumn{3}{c||}{$\boldsymbol{E_i},\, i=6,7,8$}                                                                                                          & \multicolumn{3}{c|}{$\boldsymbol{H_3}$}                                                     \\ \hline
$a$  & \multicolumn{3}{c||}{$\sigma_3$}                                                                                                                &\multicolumn{3}{c||}{$\sigma_4$}                                                                                                                 & \multicolumn{3}{c|}{$\sigma_2$} \\ \hline
$s'$    & $\sigma_1$  & $\sigma_3$ & $\Sigma_{2,n}$                                                                         & $\sigma_1$ & $\sigma_4$ & $\Sigma_{2,n}$                                                                      & $\sigma_1$  & $\sigma_2$ & $\sigma_3$                                \\ \hline
$\tilde{s}$ & $\sigma_3^2\sigma_1^3\sigma_3^2$ & $\sigma_3$ & $\,\,\Sigma_{3,2}^{(2)}\Sigma_{2,n}\Sigma_{n,3}^{(2)}\,$  &$\,\sigma_4^2\sigma_1^3\sigma_4^2\,$ & $\,\sigma_4\,$ & $\,\,\Sigma_{4,1}^{(2)}\Sigma_{2,n}\Sigma_{n,4}^{(2)}\,$ & $\sigma_2^2\sigma_1\sigma_2^2$ & $\,\sigma_2\,$ & $\sigma_2^2\sigma_1^3\sigma_2^2$\\ \hline
\end{tabular}

\medskip

\begin{tabular}{|c||c|c|c||c|c|c|c||c|c|}
\hline
$A$  & \multicolumn{3}{c||}{$\boldsymbol{H_4}$}                                                                                                                                                                                                                        & \multicolumn{4}{c||}{$\boldsymbol{F_4}$}                                                                                                & \multicolumn{2}{c|}{$\boldsymbol{I_{2m}}$}\\ \hline
$a$   & \multicolumn{3}{c||}{$\sigma_2$}                                                                                                                                                               & \multicolumn{4}{c||}{$\sigma_2$}                                                                                           & \multicolumn{2}{c|}{$\sigma_2$}\\ \hline
$s'$  & $\sigma_1$  & $\sigma_2$ & $\Sigma_{2,4}$ & $\sigma_1$ & $\sigma_2$ & $\sigma_3$  & $\sigma_4$  & $\sigma_1$  & $\sigma_2$ \\ \hline
$\tilde{s}$  & $\,\,\sigma_2^2\sigma_1\sigma_2^2\,\,$ & $\,\,\sigma_2\,$ & $\,\,\Sigma^{(2)}_{2,4}\Sigma_{4,2}\sigma_2^2\,\,$ & $\,\,\sigma_2^2\sigma_1^3\sigma_2^2\,\,$ & $\,\sigma_2\,$ & $\,\sigma_2^2\sigma_3\sigma_2^2\,$ & $\,\,\sigma_2^2\sigma_3^2\sigma_4\sigma_2^2\,\,$  & $\,\sigma_2^2\sigma_1\sigma_2^2\,$ & $\,\sigma_2$ \\ \hline
\end{tabular}

\caption{Every $\tilde{s}\in A$ projects to a generator $s'$ of the  associated Coxeter group of $A$, $W_A$. Then, with the words of the form $\tilde{s}$ we can construct, for every $p\in W_A$, an element of $A$ that projects to $p$ and whose left normal form has as first and last factor the atom $a$.}\label{tablaprueba}
\end{table}

\begin{figure}
\centering
\begin{subfigure}{0.3\linewidth}
  \centering
   \includegraphics[]{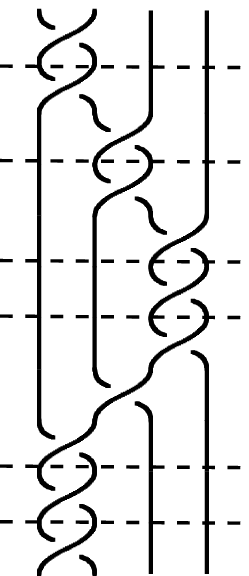}
  \caption{$ a=\sigma_1$}
\end{subfigure}%
\begin{subfigure}{0.3\linewidth}
  \centering
   \includegraphics[]{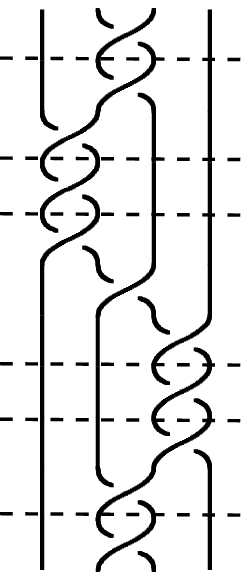}
  \caption{$ a=\sigma_2$}
\end{subfigure}%
\begin{subfigure}{0.3\linewidth}
  \centering
   \includegraphics[]{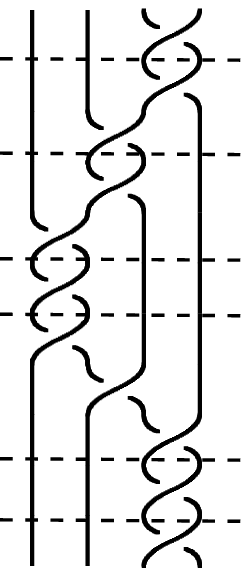}
  \caption{$a=\sigma_3$}
\end{subfigure}%
\caption{ Examples of elements $\tilde{s}$ of $A_3$ that projects to a $4$-cycle whose  left normal form have as first and last factor of the atom $a$. The dashed lines separate the factors of the normal form.}
\label{examples}
\end{figure}

\medskip
In the case where $A$ is the braid group with $n+1$ strands, $A=A_{n}$, the associated Coxeter group is the symmetric group. Then, to generate an arbitrary permutation we need a $(n+1)$-cycle and a transposition of two adjacent element in that $n$-cycle. We will chose as transposition $(1\,2),$ which will be represented in $S$ as $\sigma_1$. As $(n+1)$-cycle, we will take either $(1\, 2\, \dots\, n+1)$ or $(n+1\, n\, \dots \, 1)$, which are represented respectively by $\sigma_1\cdots\sigma_n$ and $\sigma_n\cdots\sigma_1$.  Examples of suitable preimages $\tilde{s}$ associates to the elements in $S$ are provided in Figure~\ref{examples}.

\medskip

Now notice that in the case $A=B_n$ we choose $a=\sigma_2$ and take $S'=\{\sigma_1, \sigma_2, \sigma_2\sigma_3\cdots \sigma_n\}$ because $s'=\sigma_2$ and $s'=\sigma_2\sigma_3\cdots \sigma_n$ generate a symmetric group with $n$ elements. Then we proceed with these generators as in braid groups. A similar reasoning applies also for $H_3,\,H_4,\,D_n\,E_6,\, E_7$ and $E_8$. The other elements of the form $\tilde{s}$ are computed in an explicit way. 

\medskip
Finally, notice that the length of $w_p$ is also bounded because the Coxeter group $W_A$ is finite by definition. This completes the proof. 

\hfill $\blacksquare$

%Next, let us study $A=B_n$. In this case we choose $a=\sigma_2$ and take $S'=\{\sigma_1, \sigma_2, \sigma_2\sigma_3\cdots \sigma_n\}$. Observe that $s=\sigma_2$ or $s=\sigma_2\sigma_3\cdots \sigma_n$ generate a braid group with $n$ strands. Hence proceed as described before for braid groups. A suitable preimage for
%$s=\sigma_1$ is $\tilde{s}=\sigma_2^2\sigma_1\sigma_2^2$. The same reasoning applies to $A = H_3$ and to $A=H_4$ by taking respectively $S'=\{\sigma_1,\sigma_2,\sigma_2\sigma_3\}$ and $S'=\{\sigma_1,\sigma_2,\sigma_2\sigma_3\sigma_4\}$. Analogously, for $A=I_{2m}$ consider $S'=\{\sigma_1,\sigma_2\}$.
%
%\medskip
%Now consider $A=D_n$, choose $a=\sigma_3$ and take $S'=\{\sigma_1, \sigma_3, \sigma_2\sigma_3\cdots \sigma_n\}$. As before, for $s=\sigma_3$ or $s=\sigma_2\sigma_3\cdots \sigma_n$ just proceed as for braid groups above, because they generate a braid group with $n$ strands. For $s=\sigma_1$, a suitable preimage is $\tilde{s}=\sigma_3^2\sigma_1^3\sigma_3^2$. This procedure also applies to $E_6$, $E_7$ and $E_8$ by choosing $a=\sigma_4$ and taking respectively $S'=\{\sigma_1,\sigma_4,\sigma_2\cdots\sigma_6\}$, $S'=\{\sigma_1,\sigma_4,\sigma_2\cdots\sigma_7\}$ and $S'=\{\sigma_1,\sigma_4,\sigma_2\cdots\sigma_8\}$.
%\medskip
%
%The last group is $A=F_4$, for which we consider $S=\{ \}$.

\medskip

\begin{theorem}
Let $A'$ be the commutator subgroup of $I_{2(2m+1)}$,  $A_n,\, D_n,\, E_n$ or $H_n$ equipped with a finite generator system $S$. Define $B(R):=B_{\Gamma(G,S)}(1,R)$ . Then there are constants $\varepsilon , R_0 > 0$ depending on $S$, such that for all $R > R_0$,

$$\dfrac{|\lox{G,\ca} \cap B (R) |}{|B (R)|} > \varepsilon $$

\end{theorem}

\emph{Proof.} Let $(A,\Sigma)$ be an Artin-Tits system and denote~$A'$  the commutator subgroup of~$A$. Consider the kernel~$K$ (of finite generation) of the homomorphism $e: A\longrightarrow \mathbb{Z}$ such that $e(\sigma_i)=1, \forall \sigma_i\in \Sigma$.  For all~$A$ we have that $A'\subseteq K$. However, if we set $A= I_{2(2m+1)}$, $A_n,\, D_n,\, E_n,\,H_n$, we have that $A/A'=A_{\text{ab}}=\mathbb{Z}$ \citep[Proposition 1]{Bellingeri2008}. Hence, $A'=K$, i.e, $A'$ is equal to the subgroup of elements with exponent sum equal to zero. To get the explicit generators of $A'$, notice that an element of~$A'$ in left normal form is written as $x=\Delta^{-k} x_1\cdots x_r$, where $k,r\geq 0$. If $e(\Delta)=p$, then $e(x_1)+\cdots + e(x_r)=p\cdot k$, because $e(x)=0$. Thus, we can write $$x= \Delta^{-1} a_1 \cdot \Delta^{-1} a_2 \cdots \Delta^{-1} a_k, \quad e(a_i)=p,\, a_i\in \mathcal{P}, \forall i=1,\dots, k.$$ This means that we can choose $S=\{\Delta^{-1}a \,|\, a\in \mathcal{P},\, e(a)=p\}$, which is finite. 

\medskip
Define $w_g$ to be the left normal form of an element $g\in A'$ and let $a$ be any atom such that $\varphi(w_g)\so a$. We claim that for each $g\in A'$ we can find an $x\in A'$ such that $g\cdot x$ is rigid and its left normal form contains the subword $w_a$. We will follow the same scheme as in Theorem~\ref{pure}, saying that the normal form of $g\cdot x$ has to be

\[ \Delta^{-2h}\cdot \underbrace{\overbrace{\boldsymbol{\iota(g)}}^{s\not\po}\cdots \overbrace{\boldsymbol{\varphi(g)}}^{\so a}}_{w_g} \cdot  \underbrace{\boldsymbol{a}\cdots \boldsymbol{a}}_{w_a} \cdot \underbrace{\boldsymbol{a}\cdots \boldsymbol{a}}_{w_c} \cdot \underbrace{\boldsymbol{a}\cdots (\boldsymbol{\Delta s^{-1}})}_{w_r}\]

\medskip

Let $d=-e(w_g\cdot w_a \cdot w_r)$ and recall that $e(\Delta)=p$. If we find a word $w_c$ and $h\in\mathbb{Z}$ such that $e(\Delta^{-2h}w_c)=d$, for some $h\in\mathbb{Z}$, and such that $\iota(w_c)=\varphi(w_c)=a$, then we can define $$x=\Delta^{-2h}\cdot w_a \cdot w_c \cdot w_r.$$ Hence, the left normal form of $g\cdot x$ would be $\Delta^{-2h} \cdot w_g \cdot w_a \cdot w_c \cdot w_r$, which would be rigid and contains $w_a$, as we want. If $d	\geq 0$, then $h=0$ and $w_c=a^d$ do the job. Otherwise, we choose $h=d$ and $w_c=a^{d(1+2p)}$. Finally, notice that $h$ and the length of $w_c$ depends on $d$, which depends on the length of $w_a$ and $w_r$. But we already know that $w_a$ and $w_r$ are bounded as in Theorem~\ref{pure}. Thus, $x$ is bounded and Lemma~\ref{tecnique} applies. \hfill $\blacksquare$

\medskip

\begin{remark}
In this article, we only look at the commutator subgroups of $I_{2(2m+1)}$, $A_n,\, D_n,\, E_n$ and $H_n$ because these are the only cases where the commutator subgroup is well-understood and finitely generated \citep{Orevkov2012}.

Notice that $I'_{2(2m)}$ is infinitely generated. On the other hand, $B'_3$ and $F'_4$ are finitely generated, but the question of their finite presentation is still open. $B_n\, (n>3)$ is finitely generated and finitely presented but does not fit in the scheme of the proof above. In fact, for $A=I_{2(2m)},\,B_n,\,F_4$ we have $A/A'=\mathbb{Z}^2$ \citep[Proposition 1]{Bellingeri2008}. However, we conjecture that these group have also a positive proportion of loxodromically acting elements. 

\end{remark}

\medskip

{\bf Note.} While putting the finishing touches on this paper, we learned about the existence of the preprint \citep{Yang2016}, which prove a more general result regarding the genericity of contracting elements of groups. The proofs and techniques explained in this paper, which mainly use Garside theory, have been done independently and simultaneously with Yang's article.  

\medskip

{\bf Acknowledgements. } This research was supported by a PhD contract founded by Universit\'{e} Rennes 1, Spanish Projects MTM2013-44233-P, MTM2016-76453-C2-1-P and FEDER and French-Spanish Mobility programme "M\'{e}rim\'{e}e 2015".

\bibliography{loxodromic}

\bigskip\bigskip{\footnotesize%
 \textit{ Mar\'{i}a Cumplido, UFR Math\'ematiques, Universit\'e de Rennes 1, France, and Departamento de \'Algebra, Universidad de Sevilla, Spain} \par
  \textit{E-mail address:} \texttt{\href{mailto:maria.cumplidocabello@univ-rennes1.fr}{maria.cumplidocabello@univ-rennes1.fr}}

\end{document}